\documentclass{article}

\usepackage[utf8]{inputenc}
\usepackage[italian,english]{babel}

\usepackage{geometry}	
\usepackage{fancyhdr}	
\usepackage{titlesec}		
\usepackage{indentfirst}  	

\usepackage[linktoc=all]{hyperref}	
\usepackage{tikz}				

\usepackage{amsmath}
\usepackage{amssymb}	
\usepackage{mathrsfs}	
\usepackage{amsthm}	

\theoremstyle{plain}
\newtheorem{mythm}{Theorem}[section]
\newtheorem{myprop}[mythm]{Proposition}
\newtheorem{mylemma}[mythm]{Lemma}
\newtheorem{mycor}[mythm]{Corollary}
\newtheorem{mydef}[mythm]{Definition}

\newtheorem{introthm}{Theorem}

\theoremstyle{remark}

\newcommand{\bN}{\mathbb N}

\newcommand{\bZ}{\mathbb Z}
\newcommand{\cM}{\mathcal M}
\newcommand{\fD}{\mathfrak D}

\newcommand{\sgr}{\le}

\newcommand{\Rar}{\Rightarrow}
\newcommand{\rar}{\rightarrow}

\newcommand{\ol}[1]{\overline{#1}}

\newcommand{\oc}[1]{\widehat{#1}}

\newcommand{\gen}[1]{\langle #1 \rangle}

\newcommand{\abs}[1]{|#1|}
\newcommand{\norma}[1]{\Vert#1\Vert}

\newcommand{\im}[1]{\mathrm{im}(#1)}
\newcommand{\fix}[1]{\mathrm{Fix}(#1)}

\newcommand{\rank}[1]{\mathrm{rank}(#1)}
\newcommand{\edges}[1]{\norma{#1}_\textrm{e}}

\newcommand{\cov}[1]{\mathrm{cov}(#1)}
\newcommand{\core}[1]{\mathrm{core}(#1)}
\newcommand{\bcore}[1]{\mathrm{core}_*(#1)}

\newcommand{\ffg}[1]{\mathrm{ffs}(#1)}
\newcommand{\fine}{\text{fine}}

\title{An algorithm to recognize echelon subgroups of a free group}
\author{
Dario Ascari\\
{\small \textit{Department of Mathematics, University of the Basque Country,}}\\
{\small \textit{Barrio Sarriena, Leioa, 48940, Spain}}\\
{\small e-mail: \texttt{ascari.maths@gmail.com}}\\
\\
}
\date{}

\begin{document}

\maketitle

\begin{abstract}
We provide an algorithm that, given a finite set of generators for a subgroup $H$ of a finitely generated free group $F$, determines whether $H$ is echelon or not and, in case of affirmative answer, also computes a basis with respect to which $H$ is in echelon form. This answers to a question of A. Rosenmann. We also prove, by means of a counterexample, that intersection of two echelon subgroups needs not to be echelon, answering to another question of A. Rosenmann.
\end{abstract}

\begin{center}
\small \textit{Keywords:} Free Groups, Free Factors, Whitehead's Algorithm\\
\small \textit{2020 Mathematics subject classification:} 20E05, 20E07 (20F65)
\end{center}

\section{Introduction}

Fixed-point subgroups of automorphisms of a free group have been widely studied over the years. In \cite{Ger87} Gersten showed that $\fix{\varphi}$ is finitely generated whenever $\varphi:F_n\rar F_n$ is an automorphism of a free group $F_n$ of rank $n$. The introduction of train-track graph representatives, due to Bestvina and Handel in \cite{BH92}, opened the way to a vast amount of new results in this direction, including a proof to Scott's conjecture, stating that $\rank{\fix{\varphi}}\le n$ for every automorphism $\varphi$ of $F_n$. It turns out that fixed-point subgroups behave in a very controlled manner with respect to intersections: in \cite{DV96} Dicks and Ventura introduce the notion of \textbf{inert} subgroup, i.e. a subgroup $H\sgr F_n$ such that $\rank{H\cap K}\le\rank{K}$ for every subgroup $K\sgr F_n$, and they show that the subgroup $\fix{\varphi}$ is inert, generalizing the result of Bestvina and Handel. An algorithm that given the automorphism $\varphi$ computes $\fix{\varphi}$ was found later, see \cite{BM12}. However, the converse problem, of determining whether a given subgroup $H\sgr F_n$ is the fixed-point subgroup of some automorphism, remains open. Some partial results can be found in the literature, based on the train track maps of Bestvina and Handel, but they only go in one direction, showing that if $H$ is a fixed-point subgroup then $H$ has to be of a certain form, see \cite{MV04}.

Recently, Rosenmann introduced the notion of echelon subgroup of a free group, in analogy with the concept of matrix in echelon form in linear algebra: a subgroup $H\sgr F_n$ is called \textit{echelon} if there is a basis $b_1,...,b_n$ for $F_n$ such that $\rank{H\cap\gen{b_1,...,b_i}}\le\rank{H\cap\gen{b_1,...,b_{i-1}}}+1$ for $i=1,...,n$ (see Definitions \ref{def:echelonbasis} and \ref{def:echelon}). It follows immediately from the normal form of Martino and Ventura \cite{MV04} that fixed-point subgroups of automorphisms of $F_n$ are echelon. In \cite{Ros13} Rosenmann shows that echelon subgroups are inert, i.e. they enjoy the same property of controlled behaviour with respect to intersections as fixed-point subgroups do. This places the property of being echelon as an intermediate property in the chain of implications
\begin{equation}\label{chain1}
\text{fixed-point subgroup of an automorphism} \Rar \text{echelon} \Rar \text{inert}
\end{equation}
Both the converse implications are false in general: being inert doesn't imply being echelon (as pointed out by Rosenmann in \cite{Ros13}), and being echelon doesn't imply being a fixed-point subgroup of an automorphism (since fixed-point subgroups of automorphisms are root-closed, while echelon subgroups need not to be).

In his paper \cite{Ros13} Rosenmann asks whether there is an algorithm that determines whether a given subgroup is echelon or not. We here give an affirmative answer to that question:

\begin{introthm}[Theorem \ref{algechelon}]\label{introalgechelon}
There is an algorithm that takes as input a finite set of generators for a subgroup $H\sgr F_n$ and tells us whether it is echelon or not. In case the answer is affirmative, it also computes an echelon basis for $H$.
\end{introthm}

As we pointed out before, the question of algorithmically determining whether a subgroup is the fixed-point subgroup of some automorphism is still open. Moreover, the question of algorithmically determining whether a subgroup of $F_n$ is inert is still open too (even though Ivanov \cite{Iva18} found an algorithm to compute the strictly related \textit{Walter Neumann coefficient} of a subgroup). These two facts, together with the chain of implications (\ref{chain1}), make the result of Theorem \ref{introalgechelon} particularly intriguing.

The algorithm of Theorem \ref{introalgechelon} is recursive in nature, and it is based on a fine property of Whitehead's algorithm, that has been developed by the author in \cite{Ascari}; see Section \ref{sec:algechelon} for the details. In Section \ref{sec:intersectionechelon} we prove, by means of a counterexample, that the intersection of two echelon subgroups needs not to be echelon; this answers in the negative to another question of Rosenmann \cite{Ros13}. We also prove that the property of being echelon isn't transitive, exhibiting an example of two finitely generated subgroups $K\sgr H\sgr F_3$ such that $K$ is echelon in $H$ and $H$ is echelon in $F_3$, but $K$ isn't echelon in $F_3$.

\vspace{0.2cm}

\textbf{Echelon subgroups and the Hanna Neumann conjecture.} Finding estimates on the rank of the intersection of two finitely generated subgroups $H_1,H_2\sgr F_n$ of a free group is a classical problem in combinatorial group theory. The well-known Hanna Neumann conjecture (see \cite{Neu56}) states that $\rank{H_1\cap H_2}-1\le(\rank{H_1}-1)(\rank{H_2}-1)$ for every $H_1,H_2\sgr F_n$ finitely generated subgroups. The conjecture remained open for many years, and was finally proved in 2011 independently by Friedman \cite{Fri11} and Mineyev \cite{Min12}. Recently Abdenbi \cite{Abd22} introduced a notion of \textit{generalized echelon subgroup} of a free group, showing that every echelon subgroup is also a generalized echelon subgroup, and that every generalized echelon subgroup is inert. Rosenmann's proof of inertia was based on the study a certain family of endomorphisms of a free group, called \textit{$1$-generator endomorphisms}, whereas Abdenbi's argument makes use of tools taken from the proof of Hanna Neumann conjecture. This suggests that there might be an intrinsic relation between the notion of (generalized) echelon subgroup and the study of the rank of the intersection of finitely generated subgroups in a free group.

\vspace{0.2cm}

\textbf{Echelon subgroups and $L^2$-independence.} The notion of $L^2$-independence has recently been introduced by Antol{\'i}n and Jaikin-Zapirain \cite{AJZ22} in relation to the study of (generalized) Hanna Neumann conjecture in several families of groups; we here focus on the case of free groups. A subgroup $H\sgr F_n$ is called \textit{$L^2$-independent} if the inclusion induces an injective map of $1$-dimensional $\ell^2$-homology $H_1^{(2)}(H)\rar H_1^{(2)}(F_n)$. For a free group $F_n$, the group ring $\bZ[F_n]$ can be embedded into a division ring $\fD$, and thus the condition of $L^2$-independency can be equivalently checked on the $1$-dimensional (co)homology with coefficients in $\fD$, as we now explain.

We define a \textit{$1$-cocycle} as a map $x:F_n\rar\fD$ such that $x(gh)=x(g)+g\cdot x(h)$, and we denote with $Z^1(F_n;\fD)$ the right $\fD$-module of $1$-cocycles; we define a \textit{$1$-coboundary} as a map $x:F_n\rar\fD$ of the form $x(g)=(g-1)\cdot\xi$ for some $\xi\in\fD$, and we denote with $B^1(F_n;\fD)$ the right $\fD$-module of $1$-coboundaries; we define the right $\fD$-module $H^1(F_n;\fD):=Z^1(F_n;\fD)/B^1(F_n;\fD)$. Given a basis $b_1,...,b_n$ for $F_n$, we can define the $1$-cocycles $b^*_i:F_n\rar\fD$ for $i=1,...,n$, given by $b^*_i(b_i)=1$ and $b^*_i(b_j)=0$ for $j\not=i$. We have that $b^*_1,...,b^*_n$ is a basis for the free $\fD$-module $Z^1(F_n,\fD)$, and $b^*_1\cdot(b_1-1)+...+b^*_n\cdot(b_n-1)$ is a generator for the free $\fD$-module $B^1(F_n,\fD)$; it follows that $H^1(F_n,\fD)$ is a free $\fD$-module of dimension $n-1$. A subgroup $H\sgr F_n$ induces a map $r:Z^1(F_n;\fD)\rar Z^1(H;\fD)$ given by restriction of $1$-cocycles on $F_n$ to $H$, and thus also a map $\ol{r}:H^1(F_n;\fD)\rar H^1(H;\fD)$. We have that $H$ is $L^2$-independent in $F_n$ if and only if $\ol{r}$ is surjective, and if and only if $r$ is surjective. This allows to test $L^2$-independence of a subgroup by looking only at the map of $\fD$-modules $\ol{r}$.

Fixed a basis $b_1,...,b_n$ for $F_n$ and a basis $h_1,...,h_k$ for $H$, the map $r:Z^1(F_n;\fD)\rar Z^1(H;\fD)$ can be represented by the matrix $R\in\cM_{k\times n}(\fD)$ given by $R_j^i=b^*_i(h_j)$ for $i=1,...,n$ and $j=1,...,k$. Suppose that $H$ is echelon with respect to the ordered basis $b_1,...,b_n$ for $F_n$, and for $j=0,1,...,k$ let $i(j)\ge0$ be the minimum integer such that $\rank{H\cap\gen{b_1,...,b_{i(j)}}}\ge j$: then we have that $0=i(0)<i(1)<...<i(k)\le n$ and we can choose a basis $h_1,...,h_k$ for $H$ such that $h_1,...,h_j$ is a basis for $H\cap\gen{b_1,...,b_{i(j)}}$ for $j=1,...,k$. With respect to the bases $b_1^*,...,b_n^*$ and $h_1^*,...,h_k^*$, the matrix $R$ is in echelon form, since $R_j^{i(j)}\not=0$ and $R_j^i=0$ for $i>i(j)$; this implies that the map $r:Z^1(F_n;\fD)\rar Z^1(H;\fD)$ is surjective. This shows that echelon subgroups of $F_n$ are $L^2$-independent. In \cite{AJZ22} they prove that $L^2$-independent subgroups of $F_n$ are inert. This gives us the chain of implications
$$\text{echelon} \Rar L^2\text{-independent} \Rar \text{inert.}$$

\section*{Acknowledgements}

My heartfelt thanks goes to my supervisor Martin Bridson, for fruitful discussion about the algorithm presented in this paper. The present research work was produced thanks to the funding provided by the Engineering and Physical Sciences Research Council.

\section{Preliminaries and notations}\label{Preliminaries}

In this section, we introduce the notation about graphs and we recall the properties that we are going to use along the paper. We consider graphs as combinatorial objects, according to \cite{Sta83}; for the proofs of some of the results we also rely on \cite{KM02}.

\subsection{Graphs and paths}

Along this paper, a \textbf{graph} is a quadruple $\Gamma=(V,E,\ol{\cdot},\iota)$ consisting of a set $V=V(\Gamma)$ of \textit{vertices}, a set $E=E(\Gamma)$ of \textit{edges}, a map $\ol{\cdot}:E\rar E$ called \textit{reverse} and a map $\iota:E\rar V$ called \textit{initial endpoint}; we require that, for every edge $e\in E$, we have $\ol{e}\not=e$ and $\ol{\ol{e}}=e$. For an edge $e\in E$, we denote with $\tau(e)=\iota(\ol{e})$ the \textit{terminal endpoint} of $e$.

A \textbf{map of graphs} $f:\Gamma\rar\Gamma'$ is just a couple of maps, vertices to vertices and edges to edges, preserving the reverse and the initial endpoint. A \textbf{subgraph} of a graph $\Gamma$ is a graph $\Gamma''$ such that $V(\Gamma'')\subseteq V(\Gamma)$ and $E(\Gamma'')\subseteq E(\Gamma)$ and the reverse and endpoint maps on $\Gamma''$ are the restrictions of the maps on $\Gamma$. A \textbf{pointed graph} is a pair $(\Gamma,*)$ where $\Gamma$ is a graph and $*\in V(\Gamma)$ is a vertex.

A \textbf{path} in a graph $\Gamma$ is a sequence $\sigma=(e_1,...,e_\ell)$ of edges $e_1,...,e_\ell\in E(\Gamma)$ for some integer $\ell\ge0$, such that $\tau(e_i)=\iota(e_{i+1})$ for $i=1,...,\ell-1$. The \textbf{image} of a path $\sigma$ is the smallest subgraph of $\Gamma$ containing all the edges $e_1,...,e_\ell$. A path $(e_1,...,e_\ell)$ is \textbf{reduced} if $e_i\not=\ol{e}_{i+1}$ for every $i\in\{1,...,\ell-1\}$. A path $(e_1,...,e_\ell)$ is \textbf{closed} if $\tau(e_\ell)=\iota(e_1)$. A closed path $(e_1,...,e_\ell)$ is \textbf{cyclically reduced} if it's reduced and $e_\ell\not=\ol{e}_1$.

A graph is \textbf{connected} if for every couple of vertices $v_0,v_1$ there is a path $(e_1,...,e_\ell)$ such that $\iota(e_1)=v_0$ and $\tau(e_\ell)=v_1$. A graph is a \textbf{tree} if it doesn't contain any non-trivial closed cyclically reduced path; in a tree, for every couple of vertices $v_0,v_1$ there is a unique reduced path $(e_1,...,e_\ell)$ with $\iota(e_1)=v_0$ and $\tau(e_\ell)=v_1$.

\subsection{The geometric realization}

Let $\Gamma=(V,E,\ol{\cdot},\iota)$ be a graph. The \textbf{geometric realization} $\abs{\Gamma}$ is the $1$-dimensional CW complex that has one $0$-cell for every vertex $v\in V$ and one $1$-cell for every pair of reverse edges $\{e,\ol{e}\}\subseteq E$; the $1$-cell corresponding to the pair $\{e,\ol{e}\}$ has its two endpoints glued onto the vertices $\iota(e)$ and $\tau(e)$. Notice that, if $\Gamma$ is finite, then the cardinality of $E(\Gamma)$ is an even number, since for every edge $e$ we also have its reverse $\ol{e}$, and the natural number $\abs{E(\Gamma)}/2$ coincides with the number of $1$-cells in the geometric realization $\abs{\Gamma}$.

\begin{mydef}
Let $\Gamma$ be a finite graph. Then define $\edges{\Gamma}:=\abs{E(\Gamma)}/2\in\bN$.
\end{mydef}

A map of graphs $f:\Gamma\rar\Gamma'$ induces a cellular map $\abs{f}:\abs{\Gamma}\rar\abs{\Gamma'}$ between the geometric realizations. A subgraph $\Gamma''\subseteq\Gamma$ induces a subcomplex $\abs{\Gamma''}\subseteq\abs{\Gamma}$. For a path $\sigma=(e_1,...,e_\ell)$ with $\ell\in\bN$ we define its \textbf{geometric realization} to be the continuous map $\abs{\sigma}:[0,1]\rar\abs{\Gamma}$ sending the interval $[\frac{i-1}{\ell},\frac{i}{\ell}]$ to the $1$-cell corresponding to the couple $\{e_i,\ol{e}_i\}$ for $i=1,...,\ell$. For a graph $\Gamma$, every covering space of $\abs{\Gamma}$ can be seen as the geometric realization $\abs{q}:\abs{\Gamma'}\rar\abs{\Gamma}$ of some map of graphs $q:\Gamma'\rar\Gamma$; with an abuse of notation, we will call \textit{covering map} also the map of graphs $q:\Gamma'\rar\Gamma$ (and not only its geometric realization $\abs{q}$).

For a pointed graph $(\Gamma,*)$ we define the fundamental group $\pi_1(\Gamma,*):=\pi_1(\abs{\Gamma},*)$ to be the fundamental group of its geometric realization. If $\Gamma$ is a connected finite graph, then the fundamental group $\pi_1(\Gamma,*)$ is a free group, whose rank is finite and doesn't depend on the chosen basepoint $*\in V(\Gamma)$.

\begin{mydef}
Let $\Gamma$ be a connected finite graph. Define $\rank{\Gamma}:=\rank{\pi_1(\Gamma,*)}\in\bN$.
\end{mydef}

\subsection{Core of a graph}


Let $(\Gamma,*)$ be a connected pointed graph: define its \textbf{pointed core} $\bcore{\Gamma}$ as the subgraph given by the union of the images of all the reduced paths from the basepoint to itself. Notice that $\bcore{\Gamma}$ is connected, and the inclusion $\bcore{\Gamma}\rar\Gamma$ induces an isomorphism of fundamental groups $\pi_1(\bcore{\Gamma},*)\rar\pi_1(\Gamma,*)$. It's possible to see $\Gamma$ as the union of $\bcore{\Gamma}$ and a family of pairwise disjoint trees, each tree having exactly one vertex (and no edge) in common with $\bcore{\Gamma}$.

Let $\Gamma$ be a connected graph which is not a tree: define its \textbf{core} $\core{\Gamma}$ as the subgraph given by the union of the images of all the closed cyclically reduced paths. The graph $\core{\Gamma}$ is connected. If $(\Gamma,*)$ is a pointed graph then $\core{\Gamma}$ is contained in $\bcore{\Gamma}$ as a subgraph; if the basepoint $*$ belongs to $\core{\Gamma}$ then $\core{\Gamma}=\bcore{\Gamma}$, otherwise the inclusion is strict.

\begin{mydef}
We say that a graph $\Gamma$ is \textbf{core} if $\core{\Gamma}=\Gamma$.
\end{mydef}

\subsection{Labeled graphs}

Let $F_n$ be a finitely generated free group with basis $a_1,...,a_n$. We denote with $\ol{a}_i=a_i^{-1}$ the inverse $a_i$. We denote with $\Delta_n$ the graph with one vertex $V=\{*\}$ and edges $E=\{a_1,\ol{a}_1,...,a_n,\ol{a}_n\}$. The fundamental group $\pi_1(\Delta_n,*)$ will be identified with $F_n$: for the path $\sigma_i=(a_i)$, the homotopy class of its geometric realization $\abs{\sigma_i}$ corresponds to the element $a_i\in F_n$.

A \textbf{labeled graph} is a graph $\Gamma$ together with a map of graphs $f:\Gamma\rar \Delta_n$. This means that every edge of $\Gamma$ is equipped with a label in $\{a_1,\ol{a}_1,...,a_n,\ol{a}_n\}$, according to which edge of $\Delta_n$ it is mapped to, reverse edges having inverse labels; the map $f:\Gamma\rar \Delta_n$ is called \textbf{labeling map} for $\Gamma$. Let $\Gamma_0,\Gamma_1$ be labeled graphs with labeling maps $f_0,f_1$ respectively: a map of graphs $h:\Gamma_0\rar \Gamma_1$ is called \textbf{label-preserving} if $f_1\circ h=f_0$.


By the theory of covering spaces, there is a bijection between covering spaces of $\Delta_n$ and subgroups of $F_n$. Given a subgroup $H\sgr F_n$ there is a unique (up to isomorphism) connected pointed covering space $p:(\Gamma,*)\rar(\Delta_n,*)$ such that $p_*:\pi_1(\Gamma,*)\rar(\Delta_n,*)$ induces an isomorphism between $\pi_1(\Gamma,*)$ and the subgroup $H$; notice that $\Gamma$ has a natural structure of labeled graph, with labeling map given by the covering space map $p$. For a subgroup $H\sgr F_n$ we denote $\cov{H}:=\Gamma$.

We define the \textbf{pointed core} of a subgroup $H\sgr F_n$ to be the labeled graph $\bcore{H}:=\bcore{\cov{H}}$: this comes naturally equipped with a basepoint, namely the basepoint $*$ of the covering space $(\cov{H},*)$. The map $p_*:\pi_1(\bcore{H},*)\rar\pi_1(\Delta_n,*)$ induces an isomorphism between $\pi_1(\bcore{H},*)$ and $H$.

\subsection{Folded graphs}

\begin{mydef}
Let $\Gamma$ be a labeled graph with labeling map $f:\Gamma\rar\Delta_n$. We say that $\Gamma$ is \textbf{folded} if for every pair of distinct edges $e,e'\in E(\Gamma)$ with $\iota(e)=\iota(e')$ we have $f(e)\not=f(e')$.
\end{mydef}

This means that for every vertex and for every label, there is at most one edge going out of that vertex and with that label; see also Definition 2.3 of \cite{KM02}. A labeled graph $\Gamma$ with labeling map $f:\Gamma\rar\Delta_n$ is folded if and only if for every reduced path $\sigma=(e_1,...,e_\ell)$ the word $f(e_1)...f(e_\ell)$ is reduced. Notice that every covering space of $\Delta_n$ is a folded graph, since covering maps are locally injective, and in particular $\bcore{H}$ is folded for every subgroup $H\sgr F_n$.

\begin{mylemma}\label{pi1injective}
Let $(\Gamma,*)$ be a pointed folded graph. Let $f:\Gamma\rar\Delta_n$ be the labeling map. Then the induced homomorphism $f_*:\pi_1(\Gamma,*)\rar F_n$ is injective.
\end{mylemma}
\begin{proof}
See Proposition 5.3 of \cite{Sta83}.
\end{proof}

\begin{myprop}
Let $(\Gamma,*)$ be a pointed folded graph. Let $f:\Gamma\rar\Delta_n$ be the labeling map and let $H=f_*(\pi_1(\Gamma,*))\sgr F_n$. Then we have $(\bcore{\Gamma},*)=(\bcore{H},*)$.
\end{myprop}
\begin{proof}
Let $\Gamma_0=\Gamma$. Suppose we are given a vertex $v\in V(\Gamma_0)$ and a label $a\in\{a_1,\ol{a}_1,...,a_n,\ol{a}_n\}$ such that $\Gamma_0$ contains no edge going out of $v$ and with label $a$ (i.e. there is no edge $e$ with $\iota(e)=v$ and $f(e)=a$). Then we add to $\Gamma_0$ a new vertex $u_{v,a}$ and a new edge $e_{v,a}$ (and its reverse $\ol{e}_{v,a}$), setting $\iota(e_{v,a})=v$ and $\tau(e_{v,a})=u_{v,a}$ and $f(e_{v,a})=a$. We perform this operation for all such couples $(v,a)$ at the same time, and we define $\Gamma_1$ to be the labeled graph that we get as a result.

We reiterate the construction on the graph $\Gamma_1$, and we thus define an infinite sequence of labeled graphs $\Gamma_0\subseteq\Gamma_1\subseteq\Gamma_2,...$. We set $\Omega=\Gamma_0\cup\Gamma_1\cup...$ to be the limit of the sequence.

Given a vertex $v\in V(\Gamma_n)$ and a label $a\in\{a_1,\ol{a}_1,...,a_n,\ol{a}_n\}$, we have that $\Gamma_{n+1}$ contains exactly one edge going out of $v$ and with label $a$. Moreover, for $m\ge n+2$, we aren't adding any edge adjacent to $v$ to $\Gamma_m$. It follows that, for every vertex $v\in V(\Omega)$ and for every label $a\in\{a_1,\ol{a}_1,...,a_n,\ol{a}_n\}$, there is exactly one edge of $\Omega$ going out of $v$ and with label $a$. This proves that the labeling map $f:\Omega\rar\Delta_n$ is a covering space.

It's easy to prove by induction on $n\in\bN$ that $\bcore{\Gamma_n}=\bcore{\Gamma}$. Since each path in $\Omega$ has to be contained in $\Gamma_n$ for some $n\in\bN$, we obtain that $\bcore{\Omega}=\bcore{\Gamma}$. It follows that $\pi_1(\Omega,*)=\pi_1(\Gamma,*)$, and thus $\Omega=\cov{H}$ and $\bcore{\Omega}=\bcore{H}$. The conclusion follows.
\end{proof}

This means that there is a bijection between the set of subgroups $H\sgr F_n$, and the set of pointed folded graphs $(\Gamma,*)$ with $\bcore{\Gamma}=\Gamma$. To a subgroup $H\sgr F_n$, we associate the pointed folded graph $(\Gamma,*)=(\bcore{H},*)$ which satisfies $\bcore{\Gamma}=\Gamma$. To a pointed folded graph $(\Gamma,*)$ with $\bcore{\Gamma}=\Gamma$, we associate the subgroup $H=f_*(\pi_1(\Gamma,*))\sgr F_n$ where $f:\Gamma\rar\Delta_n$ is the labeling map.

\subsection{Finitely generated subgroups and intersections}

The following Proposition \ref{fingen-finite} shows that finitely generated subgroups correspond to finite graphs.

\begin{myprop}\label{fingen-finite}
Let $H\sgr F_n$ be a subgroup. Then $H$ is finitely generated if and only if $\bcore{H}$ is finite. Moreover, given a finite set of generators for $H$, it's possible to algorithmically construct the finite pointed folded graph $\bcore{H}$.
\end{myprop}
\begin{proof}
See Lemma 5.4 in \cite{KM02} and Proposition 7.1 of \cite{KM02}. The algorithm is based on Stallings folding operations, and for a detailed description of how the algorithm works the reader can refer to Algorithm 5.4 in \cite{Sta83} or to Proposition 3.8 in \cite{KM02}.
\end{proof}

It's also possible to compute algorithmically the pointed core of the intersection of two given finitely generated subgroups of $F_n$.

\begin{myprop}\label{subgroups-intersection}
Let $H_1,H_2\sgr F_n$ be finitely generated subgroups. Then $H_1\cap H_2$ is finitely generated. Moreover, there is an algorithm that, given finite sets of generators for $H_1$ and $H_2$, computes the finite graph $\bcore{H_1\cap H_2}$.
\end{myprop}
\begin{proof}
See Theorem 5.5 and Theorem 5.6 of \cite{Sta83}.
\end{proof}

Since we will need to explicitly apply the algorithmic construction of the above Proposition \ref{subgroups-intersection}, we outline it here. Given two finitely generated subgroups $H_1,H_2\sgr F_n$, consider the two pointed graphs $(\bcore{H_1},*)$ and $(\bcore{H_2},*)$ with labeling maps $f_1:\bcore{H_1}\rar\Delta_n$ and $f_2:\bcore{H_2}\rar\Delta_n$. Consider the pull-back $\Omega=\bcore{H_1}\times_{\Delta_n}\bcore{H_2}$ as defined in Section 1.3 of \cite{Sta83}. This means that $\Omega$ has vertices $V(\Omega)=V(\bcore{H_1})\times V(\bcore{H_1})$ and edges $E(\Omega)=\{(e_1,e_2)\in E(\bcore{H_1})\times E(\bcore{H_2}) : f_1(e_1)=f_2(e_2)\}$, and the two projection maps $p_1:\Omega\rar\bcore{H}$ and $p_2:\Omega\rar\bcore{H}$ satisfy $f_1\circ p_1=f_2\circ p_2$; moreover $\Omega$ has a basepoint $(*,*)$. Then we have that $(\bcore{H_1\cap H_2},*)=\bcore{\Omega,(*,*)}$.

\subsection{Core graphs and conjugacy classes of subgroups}

We define the \textbf{core} of a non-trivial subgroup $H\sgr F_n$ to be the labeled graph $\core{H}:=\core{\cov{H}}=\core{\bcore{H}}$. It follows from the above discussion that $\core{H}$ is a connected folded graph and that $\core{H}$ is finite if and only if $H$ is finitely generated. We observe that two subgroups $H_1,H_2\sgr F_n$ are conjugate if and only if they are associated to the same covering space but with different basepoints; in particular $H_1,H_2$ are conjugate if and only if they have the same core $\core{H_1}=\core{H_2}$ (see also Proposition 7.7 in \cite{KM02}).

This means that there is a bijection between the set of conjugacy classes $[H]$ of non-trivial subgroups $H\sgr F_n$, and the set of folded core graphs. To the conjugacy class $[H]$ with $H\sgr F_n$ non-trivial, we associate the folded core graph $\core{H}$. To the folded core graph $\Gamma$, we associate the conjugacy class $[H]$ of the subgroup $H=f_*(\pi_1(\Gamma,v))\sgr F_n$ where $f:\Gamma\rar\Delta_n$ is the labeling map and $v\in V(\Gamma)$ is any vertex.

Given two subgroups $H'\sgr H\sgr F_n$, the inclusion map $i:H'\rar H$ induces a label-preserving map $\oc{i}:(\cov{H'},*)\rar(\cov{H},*)$. Since the graphs are folded, this restricts to maps between the respective (pointed) core graphs, which, by an abuse of notation, we still denote with $\oc{i}:\bcore{H'}\rar\bcore{H}$ and $\oc{i}:\core{H'}\rar\core{H}$.

\section{Free factors and Whitehead's algorithm}

In this section we focus on free factors of a free group, introducing the technical tools that we will need in the subsequent Section \ref{sec:echelon}.

\subsection{Free factors}

Let $F_n$ be a finitely generated free group with basis $a_1,...,a_n$.

\begin{mydef}
A subgroup $A\sgr F_n$ is called a \textbf{free factor} if there is a basis for $A$ can be extended to a basis for $F_n$.
\end{mydef}

If a basis for $A$ can be completed to a basis for $F_n$ by adding elements $b_1,...,b_k\in F_n$, then every basis of $A$ can be completed to a basis for $F_n$ by adding the same elements $b_1,...,b_k$.

\begin{mylemma}\label{subgraph-freefactor}
Let $\Gamma$ be a graph. Let $\Gamma'$ be a subgraph with inclusion map $h:\Gamma'\rar\Gamma$, and let $v\in V(\Gamma')$ be a vertex. Then the map $h_*:\pi_1(\Gamma',v)\rar\pi_1(\Gamma,v)$ is injective and $\pi_1(\Gamma',v)$ is a free factor in $\pi_1(\Gamma,v)$.
\end{mylemma}
\begin{proof}
See Corollary 6.3 in \cite{KM02}.
\end{proof}

\begin{myprop}\label{fftranslate}
Let $H\sgr F_n$ be a finitely generated subgroup, and let $A\sgr F_n$ be a free factor. Then $H\cap A$ is a free factor in $H$.
\end{myprop}
\begin{proof}
Without loss of generality, we can assume that $A=\gen{a_1,...,a_r}\sgr F_n$ for some $1\le r\le n$. Let $\Gamma=\bcore{H}$ with basepoint $*$ and labeling map $f:\Gamma\rar\Delta_n$. We have that $\bcore{A}$ has one vertex and $r$ edges plus their $r$ reverses, labeled $a_1,...,a_n$ and $\ol{a}_1,...,\ol{a}_r$ respectively. Using the construction of Proposition \ref{subgroups-intersection} we obtain that $H\cap A=\pi_1(\Gamma',*)$, where $\Gamma'$ is the smallest subgraph of $\Gamma$ containing all the edges with label in $\{a_1,\ol{a}_1,...,a_r,\ol{a}_r\}$. From Lemma \ref{subgraph-freefactor} we deduce that $\pi_1(\Gamma',*)=H\cap A$ is a free factor of $\pi_1(\Gamma,*)=H$.
\end{proof}

The above Proposition \ref{fftranslate} has several immediate corollaries.

\begin{mycor}\label{fftower}
Let $A\sgr F_n$ be a free factor and let $B\sgr A$ be a subgroup. Then $B$ is a free factor in $A$ if and only if $B$ is a free factor in $F_n$.
\end{mycor}

\begin{mycor}\label{fftower2}
Let $A\sgr H\sgr F_n$ be subgroups and suppose that $A$ is a free factor in $F_n$. Then $A$ is a free factor in $H$ and $\rank{H}\ge\rank{A}$.
\end{mycor}

\begin{mycor}\label{ffintersection}
Let $A,B\sgr F_n$ be free factors. Then $A\cap B\sgr F_n$ is a free factor.
\end{mycor}

\subsection{The ``injective or surjective'' lemma}

The following lemma will play a key role in the proof of Theorem \ref{introalgechelon}.

\begin{mylemma}[Injective or surjective]\label{injorsurj}
Let $K\sgr H\sgr F_n$ be finitely generated subgroups. Suppose $K$ is a free factor in $H$ with $\rank{K}=\rank{H}-1$. Let $i:K\rar H$ be the inclusion map and let $\oc{i}:\core{K}\rar\core{H}$ be the induced map between their core graphs. Then $\oc{i}$ is either injective or surjective.
\end{mylemma}
\begin{proof}
Up to conjugation, we can assume that the basepoint $*$ of $\bcore{K}$ belongs to $\im{\oc{i}}$. Let $K'=\pi_1(\im{\oc{i}},*)$ and notice that by Lemma \ref{subgraph-freefactor} we have that $K'\sgr H$ is a free factor; moreover we have an inclusion $K\sgr K'$, induced by the map $\oc{i}:\bcore{K}\rar\bcore{K'}$. Suppose $\oc{i}$ isn't surjective: then we have $\rank{K'}\le\rank{H}-1=\rank{K}$. But for free factors of $H$, the inclusion $K\sgr K'$ and the inequality $\rank{K'}\le\rank{K}$ force $K'=K$. It follows that $\core{K}=\core{K'}=\im{\oc{i}}$ is a subgraph of $\core{H}$, and thus $\oc{i}$ is injective, as desired.
\end{proof}

\subsection{Whithead automorphisms and Whitehead's algorithm}

We here review the classical Whitehead's algorithm to detect free factors. The algorithm is based on finding a sequence of Whitehead automorphisms, that we can apply to the free factor without increasing the number of edges in its core graph. We will need an improved version of such algorithm, for which we refer to \cite{Ascari}. The improved version allows to obtain control not only on the core graph of the free factor, but also on the core graph of all of its subgroups. This extra property will be of fundamental importance in developing the algorithm of Theorem \ref{introalgechelon}.

Let $F_n$ be a finitely generated free group with basis $a_1,...,a_n$.

\begin{mydef}
Let $a\in\{a_1,\ol{a}_1,...,a_n,\ol{a}_n\}$ and let $A\subseteq\{a_1,\ol{a}_1,...,a_n,\ol{a}_n\}\setminus\{a,\ol a\}$. Define the \textbf{Whitehead automorphism} $\varphi=(A,a)$ as the automorphism given by $a\mapsto a$ and

\begin{center}
$\begin{cases}
a_j\mapsto a_j & \text{if } a_j,\ol{a}_j\not\in A\\
a_j\mapsto aa_j & \text{if } a_j\in A \text{ and } \ol{a}_j\not\in A\\
a_j\mapsto a_j\ol{a} & \text{if } a_j\not\in A \text{ and } \ol{a}_j\in A\\
a_j\mapsto aa_j\ol{a} & \text{if } a_j,\ol{a}_j\in A\\
\end{cases}$
\end{center}
for $a_j\not=a,\ol{a}$.
\end{mydef}

\begin{mydef}
Let $\Gamma$ be a labeled graph with labeling map $f:\Gamma\rar\Delta_n$ and let $v\in V(\Gamma)$ be a vertex. Define the \textbf{letters at $v$} to be the subset $L(v)\subseteq\{a_1,\ol{a}_1,...,a_n,\ol{a}_n\}$ of the labels of the edges coming out of $v$. This means that for $b\in\{a_1,\ol{a}_1,...,a_n,\ol{a}_n\}$ we have $b\in L(v)$ if and only if $b=f(e)$ for some edge $e\in E(\Gamma)$ with $\iota(e)=v$.
\end{mydef}

\begin{mydef}\label{fine}
Let $H\sgr F_n$ be a non-trivial finitely generated subgroup and let $\varphi=(A,a)$ be a Whitehead automorphism. We say that the action of $\varphi$ on $H$ is \textbf{$\fine$} if for each vertex $v$ of $\core{H}$, exactly one of the following configurations takes place:

(i) $L(v)\cap A=\emptyset$.

(ii) $L(v)\subseteq A$.

(iii) $a\in L(v)$ and $L(v)\subseteq A\cup\{a\}$.
\end{mydef}

The following Theorem \ref{Whitehead} contains Whitehead's classical result (which is part (i) of the theorem, see \cite{Whi36a} and \cite{Ger84}) together with an improvement (which is part (ii) of the theorem, see \cite{Ascari}).

\begin{mythm}\label{Whitehead}
Let $H\sgr F_n$ be a free factor, and suppose $\core{H}$ has more than one vertex. Then there is a Whitehead automorphism $\varphi$ such that

(i) $\edges{\core{\varphi(H)}}<\edges{\core{H}}$.

(ii) The action of $\varphi$ on $H$ is $\fine$.
\end{mythm}
\begin{proof}
See Theorems 5.6 and 5.7 in \cite{Ascari}.
\end{proof}

Part (i) of Theorem \ref{Whitehead} gives the following corollary:

\begin{mycor}
There is an algorithm that, given a finite set of generators for a subgroup $H\sgr F_n$, determines whether $H$ is a free factor.
\end{mycor}
\begin{proof}[Algorithm]
We take the subgroup $H$ and we construct $\core{H}$ (notice that this can be done algorithmically): if $\core{H}$ has just one vertex, then $H$ is a free factor and we are done. Otherwise, we take a Whitehead automorphism $\varphi$, we construct the graph $\core{\varphi(H)}$ and we check whether $\edges{\core{\varphi(H)}}<\edges{\core{H}}$. We try all the finitely many Whitehead automorphisms, and if none of them decreases the number of edges of the core graph, then by Theorem \ref{Whitehead} we obtain that $H$ isn't a free factor. If we find a Whitehead automorphism $\varphi$ that decreases the number of edges of the core graph, then we reiterate the same procedure on $\varphi(H)$ instead of $H$. Since the number of edges in the core graph strictly decreases at each step, the algorithm must terminate.
\end{proof}

We now provide several lemmas, showing how part (ii) of Theorem \ref{Whitehead} can be used to keep under control the subgroups of $H$ while running the algorithm.

\begin{mylemma}\label{precise}
Let $H\sgr F_n$ be a non-trivial finitely generated subgroup, and let $\varphi=(A,a)$ be a Whitehead automorphism such that the action of $\varphi$ on $H$ is $\fine$. If case (iii) of Definition \ref{fine} takes places for exactly $p\ge1$ vertices $v$ of $\core{H}$, then $\edges{\core{\varphi(H)}}=\edges{\core{H}}-p$. If case (iii) of Definition \ref{fine} never happens, then $\core{\varphi(H)}$ is isomorphic to $\core{H}$.
\end{mylemma}
\begin{proof}
See Lemma 5.11 in \cite{Ascari}.
\end{proof}

\begin{mylemma}\label{finesubgroups}
Let $K\sgr H\sgr F_n$ be non-trivial finitely generated subgroups, and let $\varphi=(A,a)$ be a Whitehead automorphism. If the action of $\varphi$ on $H$ is $\fine$, then the action of $\varphi$ on $K$ is $\fine$.
\end{mylemma}
\begin{proof}
See Lemma 5.12 in \cite{Ascari}.
\end{proof}

Let $K\sgr H\sgr F_n$ be finitely generated non-trivial subgroups. Let $i:K\rar H$ be the inclusion, inducing a map of graphs $\oc{i}:\core{H}\rar\core{K}$, and consider the subgraph $\im{\oc{i}}$ of $\core{H}$. For an automorphism $\varphi:F_n\rar F_n$, let $j:\varphi(K)\rar\varphi(H)$ be the inclusion, inducing a map of graphs $\oc{j}:\core{\varphi(H)}\rar\core{\varphi(K)}$, and consider the subgraph $\im{\oc{j}}$ of $\core{\varphi(H)}$.

\begin{mylemma}\label{fineimage}
Let $K\sgr H\sgr F_n$ be non-trivial finitely generated subgroups. Let $\varphi=(A,a)$ be a Whitehead automorphism such that the action of $\varphi$ on $K$ is $\fine$. Then $\edges{\im{\oc{j}}}\le\edges{\im{\oc{i}}}$.
\end{mylemma}
\begin{proof}
See Lemma 5.13 in \cite{Ascari}.
\end{proof}

\subsection{The free factor support of a subgroup}

Let $F_n$ be a finitely generated free group with basis $a_1,...,a_n$.

\begin{myprop}\label{ffg}
Let $H\sgr F_n$ be a subgroup. Then there is a unique free factor $B\sgr F_n$ containing $H$ of minimum rank. Moreover, the free factor $B$ is minimum by inclusion among the free factors containing $H$.
\end{myprop}
\begin{proof}
Let $B\sgr F_n$ be a free factor containing $H$ and such that $B$ has minimum rank among the free factors containing $H$. Let now $B'\sgr F_n$ be any other free factor containing $H$: by Corollary \ref{ffintersection} we have that $B\cap B'\sgr F_n$ is a free factor containing $H$. But by Proposition \ref{fftranslate} we have that $B\cap B'$ is a free factor in $B$ too, implying that $\rank{B\cap B'}\le\rank{B}$. But $B$ has minimum rank among the free factors containing $H$, and thus we must have $\rank{B\cap B'}=\rank{B}$, which gives $B\cap B'=B$. Thus $B$ is minimum by inclusion among the free factors containing $H$ (and in particular, it is the unique free factor of minimum rank containing $H$).
\end{proof}

\begin{mydef}\label{def:ffg}
For a subgroup $H\sgr F_n$ define the \textbf{free factor support} of $H$ to be the free factor $B=\ffg{H}$ given by Proposition \ref{ffg}.
\end{mydef}

Given a non-trivial finitely generated subgroup $H\sgr F_n$, it's possible to algorithmically compute the free factor $\ffg{H}$, as we now explain.

\begin{myprop}\label{Whitehead6}
Let $H\sgr F_n$ be a non-trivial finitely generated subgroup. Suppose that, among the labels of the edges of $\core{H}$, there are at least $\rank{\ffg{H}}+1$ different letters from $\{a_1,...,a_n\}$. Then there is a Whitehead automorphism $\varphi$ such that $\edges{\core{\varphi(H)}}<\edges{\core{H}}$.
\end{myprop}
\begin{proof}
Let $B=\ffg{H}$. By an iterated application of Theorem \ref{Whitehead}, we can find a sequence of Whitehead automorphisms $\varphi_1,...,\varphi_k$ such that the following holds:

(i) The action of $\varphi_i$ on $\varphi_{i-1}\circ...\circ\varphi_1(B)$ is fine for $i=1,...,k$.

(ii) We have $\edges{\core{\varphi_i\circ...\circ\varphi_1(B)}}<\edges{\core{\varphi_{i-1}\circ...\circ\varphi_1(B)}}$ for $i=1,...,k$.

(iii) $\core{\varphi_k\circ...\circ\varphi_1(B)}$ is a one-vertex graph with $\rank{B}$ edges.

From (i) and from Lemma \ref{finesubgroups}, we deduce that the action of $\varphi_i$ on $\varphi_{i-1}\circ...\circ\varphi_1(H)$ is fine for $i=1,...,k$. From (iii) we obtain that $\core{\varphi_k\circ...\circ\varphi_1(H)}$ contains edges with at most $\rank{B}$ different labels, and thus must be different from $\core{H}$. From (ii) and from Lemma \ref{precise}, we have that for each $i=1,...,k$ either
$$\core{\varphi_i\circ...\circ\varphi_1(H)}=\core{\varphi_{i-1}\circ...\circ\varphi_1(H)}$$
or
$$\edges{\core{\varphi_i\circ...\circ\varphi_1(H)}}<\edges{\core{\varphi_{i-1}\circ...\circ\varphi_1(H)}}$$
In particular, we can take $l\ge 1$ to be the smallest such that $\core{\varphi_l\circ...\circ\varphi_1(H)}\not=\core{\varphi_{l-1}\circ...\circ\varphi_1(H)}$. Then we have that $\core{\varphi_{l-1}\circ...\circ\varphi_1(H)}=\core{H}$ and the Whitehead automorphism $\varphi=\varphi_l$ satisfies the thesis.
\end{proof}

\begin{mythm}\label{ffgalgorithm}
There is an algorithm that, given a finite set of generators for a subgroup $H\sgr F_n$, computes a basis for $\ffg{H}$.
\end{mythm}
\begin{proof}[Algorithm]
We look for a Whitehead automorphism $\varphi$ such that $\edges{\core{\varphi(H)}}<\edges{\core{H}}$; since there is only a finite number of Whitehead automorphisms, we can to this in a bounded amount of time. If we find such $\varphi$, then we replace $H$ with $\varphi(H)$, and reiterate the process. Since the number of edges of $\core{H}$ is decreasing strictly at each step, we eventually stop. We end up with a sequence of Whitehead automorphisms $\varphi_1,...,\varphi_k$ such that $\edges{\core{\varphi_k\circ...\circ\varphi_1(H)}}$ can't be reduced by means of a Whitehead automorphism. Let $S\subseteq\{a_1,...,a_n\}$ be the set of all the letters which appear as labels of at least one edge of $\core{\varphi_k\circ...\circ\varphi_1(H)}$. Up to conjugation, we have $\varphi_k\circ...\circ\varphi_1(H)\sgr\gen{S}$ and thus $\ffg{\varphi_k\circ...\circ\varphi_1(H)}\sgr\gen{S}$. If the inclusion is strict, then $\rank{\ffg{\varphi_k\circ...\circ\varphi_1(H)}}<\abs{S}$, and we can apply Proposition \ref{Whitehead6} and find a Whitehead automorphism reducing the number of edges of $\core{\varphi_k\circ...\circ\varphi_1(H)}$, contradiction. This shows that $\ffg{\varphi_k\circ...\circ\varphi_1(H)}=\gen{S}$ and thus $\ffg{H}=\varphi_1^{-1}\circ...\circ\varphi_k^{-1}(\gen{S})$, giving an explicit basis for $\ffg{H}$.
\end{proof}

\section{Echelon subgroups of a free group}\label{sec:echelon}

In this section we introduce the notion of echelon subgroup of a free group. We prove that the property of being echelon can be recognized algorithmically, answering a question of \cite{Ros13}. We also provide two couterexamples, showing that an echelon subgroup of an echelon subgroup is not necessarily echelon, and that intersection of echelon subgroups needs not to be echelon (answering another question of \cite{Ros13}).

\subsection{Echelon subgroups}

Let $F_n$ be a free group with basis $a_1,...,a_n$. Let $H\sgr F_n$ be a finitely generated subgroup and let $r=\rank{H}$.

\begin{mydef}
A \textbf{flag} for $H$ is a chain of free factors $B_1\sgr...\sgr B_r\sgr F_n$ such that
$$\rank{H\cap B_i}=i \qquad\qquad \text{for }i=1,...,r.$$
A flag for $H$ is called \textbf{minimal} if it satisfies $\ffg{H\cap B_i}=B_i$ for $i=1,...,r$.
\end{mydef}

\begin{mydef}\label{def:echelonbasis}
An \textbf{echelon basis} for $H$ is a (ordered) basis $b_1,...,b_n$ for $F_n$ such that
$$\rank{H\cap\gen{b_1,...,b_i}}\le\rank{H\cap\gen{b_1,...,b_{i-1}}}+1 \qquad\qquad \text{for }i=1,...,n.$$
\end{mydef}

The above definition takes inspiration from the notion of matrix in echelon form, coming from linear algebra. The following lemma shows that a subgroup admits an echelon basis if and only if it admits a (minimal) flag.

\begin{mylemma}\label{echelemma1}
Let $H\sgr F_n$ be a finitely generated subgroup with $\rank{H}=r$. Then the following are equivalent:

(i) $H$ has an echelon basis.

(ii) $H$ has a flag.

(iii) $H$ has a minimal flag.
\end{mylemma}
\begin{proof}
Denote $r=\rank{H}$.

(i)$\Rar$(ii) Let $b_1,...,b_n$ be an echelon basis for $H$. By Proposition \ref{fftranslate} we have that $H\cap\gen{b_1,...,b_{j-1}}$ is a free factor in $H\cap\gen{b_1,...,b_j}$, implying that $\rank{H\cap\gen{b_1,...,b_{j-1}}}\le\rank{H\cap\gen{b_1,...,b_j}}$. Thus we have a non-decreasing sequence of integers $0\le\rank{H\cap\gen{b_1}}\le\rank{H\cap\gen{b_1,b_2}}\le...\le\rank{H\cap\gen{b_1,...,b_{n-1}}}\le r$, and since $H$ is echelon we have that at each step the sequence increases by at most one. For each $i=1,...,r$ we can thus choose an index $j(i)$ such that $\rank{H\cap\gen{b_1,...,b_{j(i)}}}=i$. Define the free factor $B_i=\gen{b_1,...,b_{j(i)}}$, and notice that $B_1\sgr...\sgr B_{r-1}\sgr B_r$ and $\rank{H\cap B_i}=i$ for $i=1,...,r$.

(ii)$\Rar$(iii) Let $B_1\sgr...\sgr B_r\sgr F_n$ be a flag for $H$. We define the free factor $B'_i=\ffg{H\cap B_i}$ for $i=1,...,r$. Since $B_i$ is a free factor containing $H\cap B_i$, we must have $B'_i\sgr B_i$, and thus also $H\cap B'_i\sgr H\cap B_i$. Since both $H$ and $B'_i$ contain $H\cap B_i$, we must have $H\cap B_i\sgr H\cap B'_i$. It follows that $H\cap B'_i=H\cap B_i$, and in particular $\ffg{H\cap B'_i}=B'_i$ and $\rank{H\cap B'_i}=i$ for $i=1,...,r$. To conclude, notice that from the definition it's immediate that $B'_1\sgr...\sgr B'_r$. Thus the free factors $B'_1,...,B'_r$ are a minimal flag for $H$.

(iii)$\Rar$(i) Let $B_1\sgr...\sgr B_r\sgr F_n$ be a minimal flag for $H$. Choose an ordered basis for $B_1$, extend it to an ordered basis for $B_2$, and so on. We obtain an ordered basis for $F_n$, and it's easy to see that this is an echelon basis for $H$.
\end{proof}

Notice that, if $B_1\sgr...\sgr B_r$ is a minimal flag for $H$, then we have $\rank{H\cap B_r}=\rank{H}$. But by Proposition \ref{fftranslate} we have that $H\cap B_r$ is a free factor in $H$, and thus this implies $H\cap B_r=H$, yielding $H\sgr B_r$ and $\ffg{H}=B_r$.

Notice that the proof of Lemma \ref{echelemma1} is constructive, thanks to Theorem \ref{ffgalgorithm}. This means that, from an algorithmic point of view, it is equivalent to have a basis with respect to which $H$ is echelon, or a (minimal) flag for $H$; in fact, from one we can compute the other.

\begin{mydef}\label{def:echelon}
We say that a finitely generated subgroup $H\sgr F_n$ is \textbf{echelon} if it admits an echelon basis.
\end{mydef}

The above Lemma \ref{echelemma1} gives the following corollary, stating that the property of being echelon is somewhat independent on the ambient free group (but not entirely, see Proposition \ref{echecounterexample2}).

\begin{mycor}\label{corollary1}
Let $H\sgr F_n$ be a finitely generated subgroup. Let $A\sgr F_n$ be a free factor containing $H$. Then $H$ is echelon in $F_n$ if and only if $H$ is echelon in $A$.
\end{mycor}

\subsection{An algorithm to recognize echelon subgroups}\label{sec:algechelon}

Let $F_n$ be a free group with basis $a_1,...,a_n$. The aim of this section is to provide an algorithm that, given a finitely generated subgroup $H\sgr F_n$, tells us whether $H$ is echelon or not, and in case of an affirmative answer also computes an echelon basis for $H$. The algorithm we provide will be recursive in nature: the key step for the recursion is given by the following Lemma \ref{echelemma2}.

\begin{mylemma}\label{echelemma2}
Let $H\sgr F_n$ be a non-trivial finitely generated subgroup. Then the following are equivalent:

(i) $H$ is echelon.

(ii) There is a free factor $A\sgr\ffg{H}$ such that $\rank{H\cap A}=\rank{H}-1$ and $H\cap A$ is echelon.
\end{mylemma}
\begin{proof}
Let $r=\rank{H}$.

(i)$\Rar$(ii). Suppose $H$ is echelon, and let $B_1\sgr...\sgr B_r\sgr F_n$ be a minimal flag for $H$, see Lemma \ref{echelemma1}. In particular we must have $H\sgr B_r$ and $\ffg{H}=B_r$. Define $A=B_{r-1}$ and we have that $A\sgr\ffg{H}$ and $\rank{H\cap A}=r-1$. We have to prove that $H\cap A$ is echelon: but we notice that the chain of free factors $B_1\sgr...\sgr B_{r-1}\sgr F_n$ is a flag for $H$, and thus we are done by Lemma \ref{echelemma1}.

(ii)$\Rar$(i). Suppose we have a free factor $A\sgr\ffg{H}$ such that $\rank{H\cap A}=r-1$ and such that $H\cap A$ is echelon. By Lemma \ref{echelemma1} we can find a flag $B_1\sgr...\sgr B_{r-1}$ such that $\rank{(H\cap A)\cap B_i}=i$ and $\ffg{(H\cap A)\cap B_i}=B_i$ for $i=1,...,r-1$. We define $B_r=\ffg{H}$ and we notice that $B_{r-1}=\ffg{H\cap A\cap B_{r-1}}\sgr\ffg{H}=B_r$, so that we have a chain of free factors $B_1\sgr...\sgr B_{r-1}\sgr B_r$. Since $A$ is a free factor, for $i=1,...,r-1$ we have $B_i=\ffg{H\cap A\cap B_i}\sgr A$ and thus $H\cap A\cap B_i=H\cap B_i$, implying that $\rank{H\cap B_i}=i$. Of course we also have $\rank{H\cap B_r}=r$. The conclusion follows by Lemma \ref{echelemma1}.
\end{proof}

We are now looking for free factors $A\sgr\ffg{H}$ such that $\rank{H\cap A}=\rank{H}-1$. In order to do this, we make use of Theorem \ref{Whitehead} and of Lemmas \ref{finesubgroups} and \ref{fineimage}, which allow us to look instead for paths inside the finite graphs $\Lambda_r(E)$, defined as follows.

\begin{mydef}
For $r,E\ge1$ integers, define the finite oriented graph $\Lambda_r(E)$ as follows.

(i) We have one vertex for each connected core folded graph $\Gamma$ with $\rank{\Gamma}=r$ and $\edges{\Gamma}\le E$.

(ii) Suppose we are given two vertices $\Gamma,\Gamma'$ of $\Lambda_r(E)$ and a Whitehead automorphism $\varphi$ such that $\edges{\Gamma'}\le\edges{\Gamma}$ and $\pi_1(\Gamma',*'),\varphi(\pi_1(\Gamma,*))$ are conjugated subgroups of $F_n$ (and this is independent on the chosen basepoints). Then we add an edge going from $\Gamma$ to $\Gamma'$, and we label that edge with the automorphism $\varphi$.
\end{mydef}

\begin{myprop}\label{echeprop}
Let $H\sgr F_n$ be a non-trivial finitely generated subgroup with $\ffg{H}=F_n$. Let $r=\rank{H}$ and let $E=\edges{\core{H}}$. Then $H$ is echelon if and only if there are vertices $\Gamma_1$ of $\Lambda_r(E)$ and $\Gamma_2,\Gamma_3$ of $\Lambda_{r-1}(E)$ such that the following conditions hold:

(i) There is a path in $\Lambda_r(E)$ from $\core{H}$ to $\Gamma_1$.

(ii) $\Gamma_2$ is a subgraph of $\core{\varphi(\pi_1(\Gamma_1,*))}$ for some Whitehead automorphism $\varphi$ (this doesn't depend on the choice of the basepoint).

(iii) There is a path in $\Lambda_{r-1}(E)$ from $\Gamma_2$ to $\Gamma_3$.

(iv) The edges of $\Gamma_3$ use only labels from a proper subset $S\subset\{a_1,...,a_n\}$.

(v) The group $\pi_1(\Gamma_3,*)$ is echelon (this doesn't depend on the choice of the basepoint).
\end{myprop}

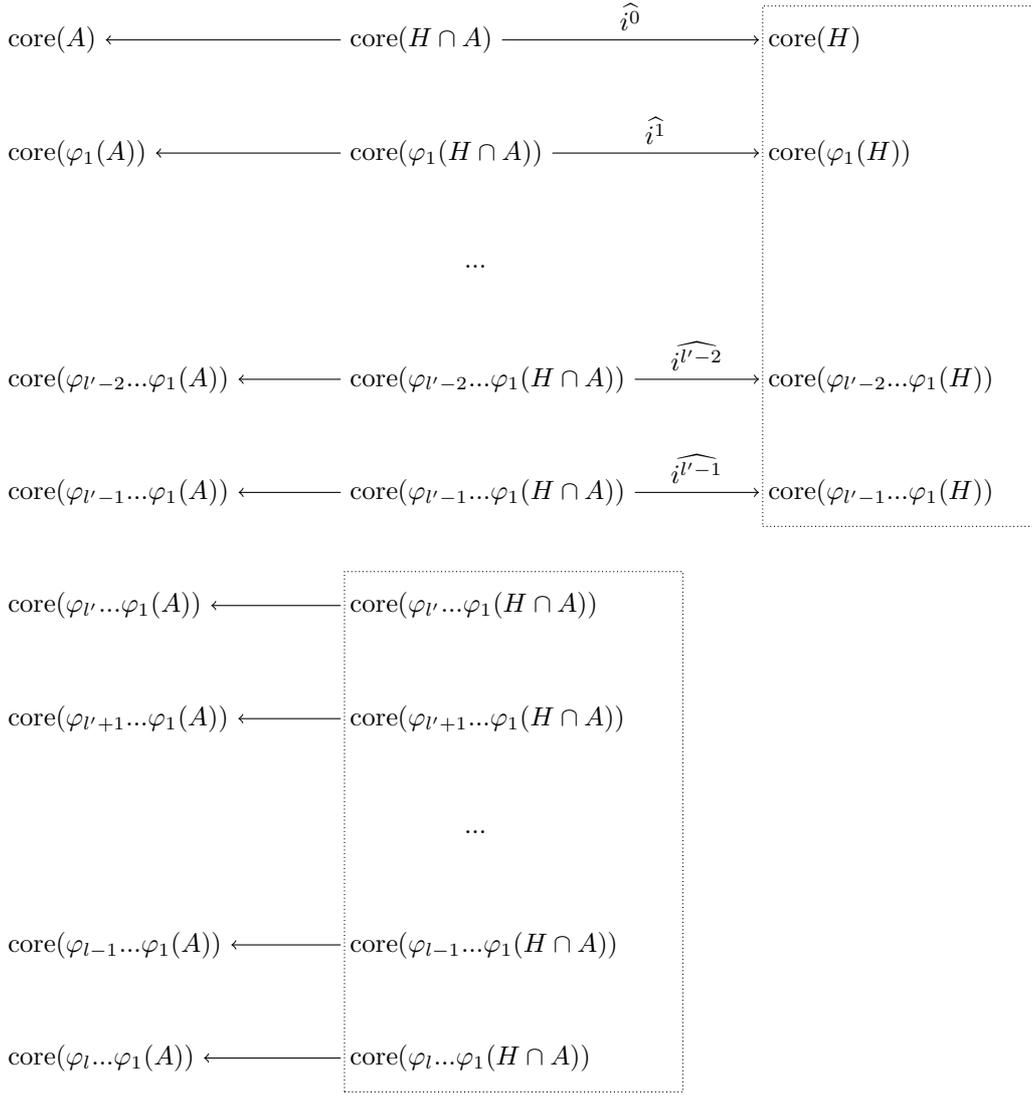
\begin{figure}[h!]
\centering
\begin{tikzpicture}[scale=1]

\draw (-3.5,0) node[right] (0a) {$\core{A}$};
\draw (1,0) node[right] (0b) {$\core{H\cap A}$};
\draw (6.5,0) node[right] (0c) {$\core{H}$}; 
\draw[->] (0b) to (0a);
\draw[->] (0b) to node[above]{$\oc{i^0}$} (0c);

\draw (-3.5,-1.5) node[right] (1a) {$\core{\varphi_1(A)}$};
\draw (1,-1.5) node[right] (1b) {$\core{\varphi_1(H\cap A)}$};
\draw (6.5,-1.5) node[right] (1c) {$\core{\varphi_1(H)}$}; 
\draw[->] (1b) to (1a);
\draw[->] (1b) to node[above]{$\oc{i^1}$} (1c);

\node[right] (2) at (2.5,-3) {...};

\draw (-3.5,-4.5) node[right] (3a) {$\core{\varphi_{l'-2}...\varphi_1(A)}$};
\draw (1,-4.5) node[right] (3b) {$\core{\varphi_{l'-2}...\varphi_1(H\cap A)}$};
\draw (6.5,-4.5) node[right] (3c) {$\core{\varphi_{l'-2}...\varphi_1(H)}$}; 
\draw[->] (3b) to (3a);
\draw[->] (3b) to node[above]{$\oc{i^{l'-2}}$} (3c);	

\draw (-3.5,-6) node[right] (4a) {$\core{\varphi_{l'-1}...\varphi_1(A)}$};
\draw (1,-6) node[right] (4b) {$\core{\varphi_{l'-1}...\varphi_1(H\cap A)}$};
\draw (6.5,-6) node[right] (4c) {$\core{\varphi_{l'-1}...\varphi_1(H)}$}; 
\draw[->] (4b) to (4a);
\draw[->] (4b) to node[above]{$\oc{i^{l'-1}}$} (4c);

\draw[densely dotted] (6.55,0.45) rectangle (10.2,-6.45);

\draw (-3.5,-7.5) node[right] (5a) {$\core{\varphi_{l'}...\varphi_1(A)}$};
\draw (1,-7.5) node[right] (5b) {$\core{\varphi_{l'}...\varphi_1(H\cap A)}$};
\draw[->] (5b) to (5a);

\draw (-3.5,-9) node[right] (6a) {$\core{\varphi_{l'+1}...\varphi_1(A)}$};
\draw (1,-9) node[right] (6b) {$\core{\varphi_{l'+1}...\varphi_1(H\cap A)}$};
\draw[->] (6b) to (6a);

\node[right] (7) at (2.5,-10.5) {...};

\draw (-3.5,-12) node[right] (8a) {$\core{\varphi_{l-1}...\varphi_1(A)}$};
\draw (1,-12) node[right] (8b) {$\core{\varphi_{l-1}...\varphi_1(H\cap A)}$};
\draw[->] (8b) to (8a);

\draw (-3.5,-13.5) node[right] (9a) {$\core{\varphi_{l}...\varphi_1(A)}$};
\draw (1,-13.5) node[right] (9b) {$\core{\varphi_{l}...\varphi_1(H\cap A)}$};
\draw[->] (9b) to (9a);

\draw[densely dotted] (1.05,-7.05) rectangle (5.5,-13.95);

\end{tikzpicture}
\caption{A diagram of how the situation in the proof of Proposition \ref{echeprop} changes as we apply the sequence of Whitehead automorphisms $\varphi_1,...,\varphi_l$. In the column on the left we can see the graph $\core{A}$, which decreases its size at each step, until we get a rose $\core{\varphi_l...\varphi_1(A)}$. In the central column we have $\core{H\cap A}$ and in the right column we have $\core{H}$. The map $\oc{i^k}$ is surjective for $k=0,...,l'-1$. The upper box describes a path in $\Lambda_r(E)$ from $\core{H}$ to $\Gamma_1$. The lower box describes a path in $\Lambda_{r-1}(E)$ from $\Gamma_2$ to $\Gamma_3$.}\label{bigdiagram}
\end{figure}

\begin{proof}\

($\Rar$) Suppose first that $H$ is echelon, and we want to find $\Gamma_1,\Gamma_2,\Gamma_3$ satisfying (i)-(v).

By Lemma \ref{echelemma2}, let $A\sgr F_n$ be a free factor such that $\rank{H\cap A}=r-1$ and $H\cap A$ is echelon. By an iterated application of Theorem \ref{Whitehead}, we find a sequence of Whitehead automorphisms $\varphi_1,...,\varphi_l$ such that for every $k=0,...,l-1$ the action of $\varphi_{k+1}$ on $\varphi_k\circ...\circ\varphi_1(A)$ is $\fine$, and such that $\core{\varphi_l\circ...\circ\varphi_1(A)}$ is a graph with one vertex. Without loss of generality, we can assume that the edges of $\core{\varphi_l\circ...\circ\varphi_1(A)}$ are labeled with the elements of $S$ for some proper subset $S\subset\{a_1,...,a_n\}$.

For $k=0,...,l$ let $i^k:\varphi_k\circ...\circ\varphi_1(H\cap A)\rar\varphi_k\circ...\circ\varphi_1(H)$ be the inclusion map, and let $\oc{i^k}:\core{\varphi_k\circ...\circ\varphi_1(H\cap A)}\rar\core{\varphi_k\circ...\circ\varphi_1(H)}$ be the induced map between core graphs, see Figure \ref{bigdiagram}. By Lemma \ref{finesubgroups} the action of $\varphi_{k+1}$ on $\varphi_k\circ...\circ\varphi_1(H\cap A)$ is $\fine$. By Lemma \ref{precise} we have $\edges{\core{\varphi_{k+1}\circ...\circ\varphi_1(H\cap A)}}\le\edges{\varphi_k\circ...\circ\varphi_1(H)}$ for $k=0,...,l-1$. By Lemma \ref{fineimage} we have $\edges{\im{\oc{i^{k+1}}}}\le\edges{\im{\oc{i^k}}}$ for $k=0,...,l-1$.

By Lemma \ref{injorsurj}, we have that $\oc{i^k}$ is either surjective or injective. Notice that $\oc{i^l}$ is not surjective: since $\ffg{H}=F_n$, we have that the edges of $\core{\varphi_l\circ...\circ\varphi_1(H)}$ make use of all the labels $\{a_1,...,a_n\}$, while $\varphi_l\circ...\circ\varphi_1(H\cap A)$ only uses labels from $S$. Let $0\le l'\le l$ be the smallest index such that $\oc{i^{l'}}$ isn't surjective. Define $\Gamma_1=\im{\oc{i^{l'-1}}}=\core{\varphi_{l'-1}\circ...\circ\varphi_1(H)}$, or $\Gamma_1=\core{H}$ if $\oc{i^0}$ is not surjective. Since $\oc{i^{l'}}$ is not surjective, it has to be injective: thus we can define $\Gamma_2=\im{\oc{i^{l'}}}=\core{\varphi_{l'}\circ...\circ\varphi_1(H\cap A)}$. Define $\Gamma_3=\core{\varphi_l\circ...\circ\varphi_1(H\cap A)}$.

For $k=0,...,l'-2$ we have $\edges{\im{\oc{i^{k+1}}}}\le\edges{\im{\oc{i^k}}}$, giving a path in $\Lambda_r(E)$ from $\im{\oc{i^0}}=\core{H}$ to $\im{\oc{i^{l'-1}}}=\Gamma_1$. This yields property (i) of the statement.

Since $\oc{i^{l'}}$ is injective, we have that $\im{\oc{i^{l'}}}=\Gamma_2$ is a subgraph of $\core{\varphi_{l'}\circ...\circ\varphi_1(H)}=\core{\varphi_{l'}(\pi_1(\Gamma_1,*)}$ and satisfies $\edges{\Gamma_2}=\edges{\im{\oc{i^{l'}}}}\le\edges{\im{\oc{i^{l'-1}}}}\le E$. This yields property (ii) of the statement.

For $k=l',...,l-1$ we have $\edges{\core{\varphi_{k+1}\circ...\circ\varphi_1(H\cap A)}}\le\edges{\varphi_k\circ...\circ\varphi_1(H)}$, giving a path in $\Lambda_{r-1}(E)$ from $\Gamma_2$ to $\Gamma_3$. This yields property (iii) of the statement.

We have an inclusion $\varphi_l\circ...\circ\varphi_1(H\cap A)\sgr\varphi_l\circ...\circ\varphi_1(A)$. Since $\core{\varphi_l\circ...\circ\varphi_1(A)}$ only uses labels from the proper subset $S\subset\{a_1,...,a_n\}$, then so does $\Gamma_3=\core{\varphi_l\circ...\circ\varphi_1(H\cap A)}$. This yields property (iv) of the statement.

The subgroup $\varphi_l\circ...\circ\varphi_1(H\cap A)$ is echelon since $H\cap A$ was by assumption. This yields property (v) of the statement.

($\Leftarrow$) Suppose we are given $\Gamma_1,\Gamma_2,\Gamma_3$ satisfying (i)-(v), and we want to prove that $H$ is echelon.

Let $B=\gen{S}$ be generated by the letters in condition (iv), so that $B$ is a proper free factor of $F_n$; take a subgroup $K_3\sgr F_n$ such that $\core{K_3}=\Gamma_3$ and $K_3\sgr B$: this can be done thanks to condition (iv). Condition (v) tells us that $K_3$ is echelon. Take also $K_1,K_2\sgr F_n$ such that $\core{K_1}=\Gamma_1$ and $\core{K_2}=\Gamma_2$.

Condition (iii) tells us that $\mu(K_3)=K_2$ for some automorphism $\mu$ of $F_n$. To be precise, we take a path in $\Lambda_{r-1}(E)$ from $\Gamma_2$ to $\Gamma_3$, we look at the Whitehead automorphisms that label the edges of that path, and we take the composition of their inverses; we also compose with a suitably chosen conjugation automorphism. Condition (ii) tells us that $\nu(K_2)$ is a free factor in $K_1$ for some automorphism $\nu$ of $F_n$. The automorphism $\nu$ is obtained as the inverse of the automorphism $\varphi$ given by condition (ii), composed with a conjugation automorphism. Condition (i) implies that $\eta(K_1)=H$ for some automorphism $\eta$ of $F_n$.

Let $A=\eta\circ\nu\circ\mu(B)$ and notice that $A$ is a proper free factor of $\ffg{H}=F_n$; in particular $H\cap A$ is a proper free factor in $H$. Let $H'=\eta\circ\nu\circ\mu(K_3)$ and notice that $H'\sgr A$ (because $K_3\sgr B$) and $H'$ is a proper free factor in $H$ (because $\nu(K_2)$ is a proper free factor in $K_1$) and $H'$ has rank $r-1$ (because $K_3$ has rank $r-1$, since $\Gamma_3$ is a vertex of $\Lambda_{r-1}(E)$). Now $H'\sgr H\cap A$ are two proper free factors of $H$ and $H'$ has rank $r-1$: this forces $H'=H\cap A$.

Thus we have found a free factor $A\sgr\ffg{H}$ such that $\rank{H\cap A}=r-1$ and $H\cap A$ is echelon. Lemma \ref{echelemma2} gives us the desired conclusion.
\end{proof}

\begin{mythm}\label{algechelon}
There is an algorithm that takes as input a finite set of generators for a subgroup $H\sgr F_n$ and tells us whether it is echelon or not. In case the answer is affirmative, it also computes an echelon basis for $H$.
\end{mythm}
\begin{proof}
We proceed by induction on $r=\rank{H}$. For $r=1$ we have that every subgroup is echelon with respect to any basis. Suppose now we have the algorithm for subgroups of rank $r-1$. Let $H\sgr F_n$ be a finitely generated subgroup of rank $\rank{H}=r$. Without loss of generality we can assume $\ffg{H}=F_n$ (otherwise we compute a basis for $\ffg{H}$ using the algorithm of Theorem \ref{ffgalgorithm}, and then we use $\ffg{H}$ as ambient group instead of $F_n$). Let $E=\edges{\core{H}}$.

We build the finite graphs $\Lambda_r(E)$ and $\Lambda_{r-1}(E)$; we take all graphs $\Gamma_1$ such that there is a path in $\Lambda_r(E)$ from $\core{H}$ to $\Gamma_1$; for each such $\Gamma_1$ and for each Whitehead automorphism $\varphi$, we take all the subgraphs $\Gamma_2$ of $\core{\varphi(\pi_1(\Gamma_1,*))}$ with at most $E$ edges; for each such $\Gamma_2$, we take all graphs $\Gamma_3$ such that there is a path in $\Lambda_{r-1}(E)$ from $\Gamma_2$ to $\Gamma_3$; finally, for each such $\Gamma_3$, we check whether it uses all the labels or not, and whether its fundamental group is echelon or not (and we are able to do this algorithmically, by inductive hypothesis, since $\rank{\Gamma_3}=r-1$). If no graph $\Gamma_3$ satisfies both the conditions, then by Proposition \ref{echeprop} we have that $H$ is not echelon, and the algorithm gives negative answer. If there is a graph $\Gamma_3$ that satisfies both the conditions, then by Proposition \ref{echeprop} we have that $H$ is echelon and the algorithm gives affirmative answer.

It remains, in case of an affirmative answer, to compute a basis with respect to which $H$ is echelon. We take the graph $\Gamma_3$ with echelon fundamental group, and we consider an echelon subgroup $K_3\sgr F_n$ with $\core{K_3}=\Gamma_3$. The inductive hypothesis on the subgroup $K_3$ allows us to compute a minimal flag $B_1\sgr...\sgr B_{r-1}\sgr F_n$ for $K_3$.

We take the proper subset $S\subset\{a_1,...,a_n\}$ of letters which appear as labels of at least one edge of $\Gamma_3$, and we call $B=\gen{S}$. We take finitely generated subgroups $K_1,K_2\sgr F_n$ with $\core{K_1}=\Gamma_1$ and $\core{K_2}=\Gamma_2$. As in the proof of Proposition \ref{echeprop}, we have $K_3\sgr B$ and we are able (by looking at the paths in $\Lambda_r(E)$ and $\Lambda_{r-1}(E)$) to compute an automorphism $\rho$ of $F_n$ such that $\rho(K_3)\sgr H$ is a proper free factor in $H$ of rank $r-1$. We finally proceed as in the proof of Lemma \ref{echelemma2}, and we consider the chain of free factors $\rho(B_1)\sgr...\sgr\rho(B_{r-1})\sgr\ffg{H}=F_n$: this is a flag for $H$, as desired.
\end{proof}

\subsection{About the intersection of echelon subgroups}\label{sec:intersectionechelon}

In \cite{Ros13} A. Rosenmann asked whether the intersection of two echelon subgroups is always echelon. We give negative answer to this question: the following proposition provides a counterexample.

\begin{figure}[h!]
\centering
\begin{tikzpicture}
\node (0) at (0,0) {$*$};

\node (1) at (-1,0) {.};
\node (2) at (-2,0) {.};
\node (3) at (-3,1) {.};
\node (4) at (-4,2) {.};
\node (5) at (-3,2) {.};
\node (6) at (-2,2) {.};
\node (7) at (-1,1) {.};

\draw[->] (0) to node[below]{$a$} (1);
\draw[->] (1) to node[below]{$a$} (2);
\draw[->] (2) to node[below]{$b$} (3);
\draw[->] (3) to node[below]{$b$} (4);
\draw[->] (4) to node[above]{$a$} (5);
\draw[->] (5) to node[above]{$a$} (6);
\draw[->] (6) to node[above]{$b$} (7);
\draw[->] (7) to node[above]{$b$} (0);

\node (8) at (1,0) {.};
\node (9) at (2,0) {.};
\node (10) at (3,1) {.};
\node (11) at (4,2) {.};
\node (12) at (3,2) {.};
\node (13) at (2,2) {.};
\node (14) at (1,1) {.};

\draw[->] (0) to node[below]{$b$} (8);
\draw[->] (8) to node[below]{$b$} (9);
\draw[->] (9) to node[below]{$a$} (10);
\draw[->] (10) to node[below]{$a$} (11);
\draw[->] (11) to node[above]{$b$} (12);
\draw[->] (12) to node[above]{$b$} (13);
\draw[->] (13) to node[above]{$a$} (14);
\draw[->] (14) to node[above]{$a$} (0);
\end{tikzpicture}
\caption{The graph $\bcore{K}$ where $K$ is generated by $a^2b^2a^2b^2$ and $b^2a^2b^2a^2$.}\label{abab}
\end{figure}
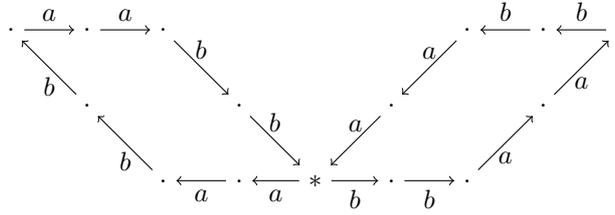

\begin{myprop}
Consider the free group $F_2=\gen{a,b}$. Then we have the following:

(i) The subgroup $H=\gen{a^2,b^2a^2b^2}$ is echelon with respect to the ordered basis $a,b$.

(ii) The subgroup $H'=\gen{b^2,a^2b^2a^2}$ is echelon with respect to the ordered basis $b,a$.

(iii) The subgroup $H\cap H'=\gen{a^2b^2a^2b^2,b^2a^2b^2a^2}$ is not echelon.
\end{myprop}
\begin{proof}
Parts (i) and (ii) are immediate. Using the algorithm of Proposition \ref{subgroups-intersection} we can see that $\bcore{H\cap H'}$ is as in Figure \ref{abab}, and thus $H\cap H'=\gen{a^2b^2a^2b^2,b^2a^2b^2a^2}$. To prove that $H\cap H'$ is not echelon, we might just use the algorithm of Theorem \ref{algechelon}, but we provide here a shorter argument. Let $K=H\cap H'$.

Suppose by contradiction that $K$ is echelon with respect to the ordered basis $p,q$ of $F_2$. Then we must have $\rank{K\cap\gen{p}}=1$ meaning that $K$ contains a power $p^\alpha$ of the primitive element $p$. The inclusion $i:\gen{p^\alpha}\rar K$ gives a map of graphs $\oc{i}:\core{\gen{p^\alpha}}\rar\core{K}$ and the image $\im{\oc{i}}$ contains at least one of the left cycle $aabbaabb$ or the right cycle $bbaabbaa$ of Figure \ref{abab}. Suppose $\core{\gen{p^\alpha}}$ contains a path $(e_1,...,e_8)$ with labels $a,a,b,b,a,a,b,b$ respectively (the other case is analogous).

By Theorem \ref{Whitehead} and Lemma \ref{finesubgroups}, there is a Whitehead automorphism $\varphi=(D,d)$, with $d\in\{a,\ol{a},b,\ol{b}\}$ and $D\subseteq\{a,\ol{a},b,\ol{b}\}\setminus\{d,\ol{d}\}$, such that the action of $\varphi$ on $\core{\gen{p^\alpha}}$ is fine (see Definition \ref{fine}). Suppose $d=a$ (the other cases are analogous): we look at the vertex $\iota(e_3)=\iota(\ol{e}_2)$ and since $\ol{a}\not\in D$ we obtain that $b\not\in D$, we look at the vertex $\iota(e_4)=\iota(\ol{e}_3)$ and we obtain that $\ol{b}\not\in D$; it follows that $D$ is empty, contradiction. The cases $d=\ol{a},b,\ol{b}$ are analogous. But then there is no Whitehead automorphism whose action on $\core{\gen{p^\alpha}}$ is fine, contradiction.
\end{proof}

%
%

The next proposition shows, by means of a counterexample, that the property of being echelon isn't transitive: there are finitely generated subgroups $K\sgr H\sgr F_n$ such that $K$ is echelon in $H$ and $H$ is echelon in $F_n$, but $K$ isn't echelon in $F_n$.

\begin{figure}[h!]
\centering
\begin{tikzpicture}
\node (0) at (0,0) {$*$};

\node (1) at (-1,0) {.};
\node (2) at (-2,0) {.};
\node (3) at (-3,0) {.};
\node (4) at (-4,0) {.};
\node (5) at (-4,1) {.};
\node (6) at (-4,2) {.};
\node (7) at (-4,3) {.};
\node (8) at (-4,4) {.};
\node (9) at (-3,4) {.};
\node (10) at (-2,4) {.};
\node (11) at (-1,4) {.};
\node (12) at (0,4) {.};

\node (13) at (1,4) {.};
\node (14) at (2,4) {.};
\node (15) at (3,4) {.};
\node (16) at (4,4) {.};
\node (17) at (5,3) {.};
\node (18) at (5,2) {.};
\node (19) at (5,1) {.};
\node (20) at (5,0) {.};
\node (21) at (4,0) {.};
\node (22) at (3,0) {.};
\node (23) at (2,0) {.};
\node (24) at (1,0) {.};

\draw[->] (0) to node[above]{$c$} (1);
\draw[->] (1) to node[above]{$c$} (2);
\draw[->] (2) to node[above]{$a$} (3);
\draw[->] (3) to node[above]{$a$} (4);
\draw[->] (4) to node[left]{$b$} (5);
\draw[->] (5) to node[left]{$b$} (6);
\draw[->] (6) to node[left]{$a$} (7);
\draw[->] (7) to node[left]{$a$} (8);
\draw[->] (8) to node[above]{$b$} (9);
\draw[->] (9) to node[above]{$b$} (10);
\draw[->] (10) to node[above]{$c$} (11);
\draw[->] (11) to node[above]{$c$} (12);
\draw[->] (12) to node[above]{$c$} (13);
\draw[->] (13) to node[above]{$c$} (14);
\draw[->] (14) to node[above]{$a$} (15);
\draw[->] (15) to node[above]{$a$} (16);
\draw[->] (16) to node[left]{$b$} (17);
\draw[->] (17) to node[right]{$b$} (18);
\draw[->] (18) to node[right]{$a$} (19);
\draw[->] (19) to node[right]{$a$} (20);
\draw[->] (20) to node[above]{$b$} (21);
\draw[->] (21) to node[above]{$b$} (22);
\draw[->] (22) to node[above]{$c$} (23);
\draw[->] (23) to node[above]{$c$} (24);
\draw[->] (24) to node[above]{$b$} (0);
\draw[->] (0) to node[right]{$a$} (12);
\end{tikzpicture}
\caption{The graph $\bcore{K}$ where $K=\gen{w\ol{a},awb}$ and $w=c^2a^2b^2a^2b^2c^2$.}\label{cc}
\end{figure}
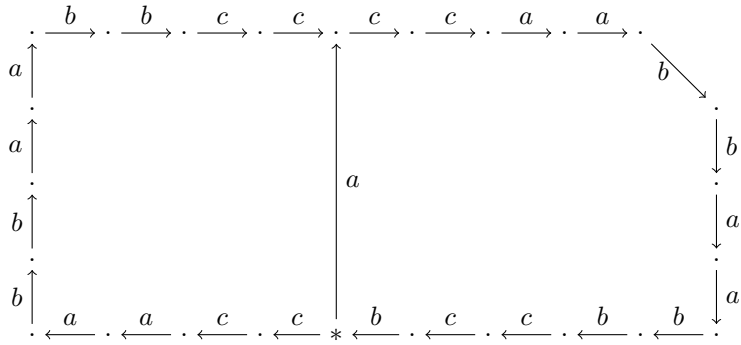

\begin{myprop}\label{echecounterexample2}
Consider the free group $F_3=\gen{a,b,c}$ and let $w=c^2a^2b^2a^2b^2c^2$. Then we have the following:

(i) The subgroup $H=\gen{a,b,w}$ is echelon in $F_3$.

(ii) The subgroup $K=\gen{w\ol{a},awb}$ is a free factor in $H$, and thus $K$ is echelon in $H$.

(iii) The subgroup $K$ isn't echelon in $F_3$.
\end{myprop}
\begin{proof}
Part (i) is immediate since $H$ is echelon with respect to the ordered basis $a,b,c$. Part (ii) is immediate too since $w\ol{a},awb$ is a basis for $K$ and $w\ol{a},awb,w$ is a basis for $H$. To prove that $K$ is not echelon in $F_3$, we may just use the algorithm of Theorem \ref{algechelon}, but we provide here a shorter argument.

Suppose by contradiction $K$ is echelon in $F_3$: then there is a proper free factor $J\sgr F_3$ such that $\rank{K\cap J}=1$. The inclusion $K\cap J\rar K$ induces a label-preserving map $\core{K\cap J}\rar\core{K}$. You can see the graph $\core{K}$ in Figure \ref{cc}, and it's easy to see that every folded core graph that has a map towards $\core{K}$ must contain a path labeled with $w=c^2a^2b^2a^2b^2c^2$. In particular $\core{K\cap J}$ has to contain a path $(e_1,...,e_{12})$ with labels $c,c,a,a,b,b,a,a,b,b,c,c$ respectively.

By Theorem \ref{Whitehead} there is a Whitehead automorphism $\varphi=(D,d)$, with $d\in\{a,\ol{a},b,\ol{b},c,\ol{c}\}$ and $D\subseteq\{a,\ol{a},b,\ol{b},c,\ol{c}\}\setminus\{d,\ol{d}\}$, such that the action of $\varphi$ on $\core{K\cap J}$ is fine (see Definition \ref{fine}). Suppose $d=a$: we look at the vertex $\iota(e_5)=\iota(\ol{e}_4)$ and since $\ol{a}\not\in D$ we obtain that $b\not\in D$, we look at $\iota(e_6)$ and we obtain that $\ol{b}\not\in D$, we look at $\iota(e_{11})$ and we obtain that $c\not\in D$, we look at $\iota(e_{12})$ and we obtain $\ol{c}\not\in D$; it follows that $D$ is empty, contradiction. The cases $d=\ol{a},b,\ol{b}$ are analogous. Suppose $d=c$: we look at $\iota(e_3)$ at we obtain that $a\not\in D$, we look at $\iota(e_4)$ and we obtain that $\ol{a}\not\in D$, we look at $\iota(e_5),\iota(e_7)$ and we obtain that $b,\ol{b}\not\in D$; it follows that $D$ is empty, contradiction. The case $d=\ol{c}$ is analogous. But then there is no Whitehead automorphism whose action on $\core{K\cap J}$ is fine, contradiction.
\end{proof}

%
%

\bibliographystyle{alpha}
\nocite{*}
\bibliography{bibliography.bib}

\end{document}